\documentclass[notitlepage,leqno,10pt]{article}

\usepackage{mathtools}

\usepackage{mathrsfs}
\usepackage{color}

\usepackage{latexsym}
\usepackage{amsmath}
\usepackage{amsfonts}
\usepackage{amssymb}
\usepackage{esint}
\usepackage{graphicx}
\usepackage{color}  

\renewcommand{\d}{\delta }
\newcommand{\D }{\Delta }

\renewcommand{\l }{\lambda }

\newcommand{\n }{\nabla }

\renewcommand{\th }{\theta }

\newcommand{\pa }{\partial}

\newcommand{\be}{\begin{equation}}
\newcommand{\ee}{\end{equation}}
\newenvironment{pf}{\noindent{\sc Proof}.\enspace}{\rule{2mm}{2mm}\medskip}
\newenvironment{pfn}{\noindent{\sc Proof}}{\rule{2mm}{2mm}\medskip}

\newcommand{\R}{\mathbb{R}}

\newcommand{\Z}{\mathbb{Z}}

\author{Sagun Chanillo\thanks{Department of Mathematics - Hill Center
		Rutgers, The State University of New Jersey, 110 Frelinghuysen Rd.
		Piscataway, NJ 08854-8019. e-mail: chanillo@math.rutgers.edu}, Andrea Malchiodi\thanks{Scuola Normale Superiore, Piazza dei Cavalieri 7, 56126 Pisa. e-mail: andrea.malchiodi@sns.it} }

\date{}

\title{Sharp bounds on the Nusselt number  in Rayleigh-B\'enard \\ convection 
 and a bilinear estimate via Carleson measures 
}

\begin{document}

\newtheorem{lem}{Lemma}[section]
\newtheorem{pro}[lem]{Proposition}
\newtheorem{thm}[lem]{Theorem}
\newtheorem{rem}[lem]{Remark}
\newtheorem{cor}[lem]{Corollary}
\newtheorem{df}[lem]{Definition}

\maketitle

\medspace \medspace

{\footnotesize
\begin{abstract}

\noindent We prove a conjecture in fluid dynamics concerning optimal bounds for 
heat transportation in the infinite Prandtl number limit. Due to a {\em maximum principle} property  
for the temperature exploited by Constantin-Doering and Otto-Seis, this amounts to proving a-priori bounds 
for horizontally-periodic solutions of a fourth-order equation in a  strip of large width. Such bounds are obtained here 
using Fourier analysis, integral representations, and a bilinear estimate due to Coifman and Meyer 
which uses the Carleson measure characterization of BMO functions by Fefferman.

\bigskip


\end{abstract}}

\section{Introduction}

 Thermal convection processes are important to understand, since they 
  appear in many applications such as engineering, meteorology, geophysics, astrophysics and oceanography. 
 Rayleigh-B\'ernard's model is a fundamental one in Convection Theory, where one considers a fluid layer between rigid plates which are heated from below and cooled from above. 
 
 We will use the convention from \cite{OS}, presenting non-dimensionalized equations: in this case the 
 relevant parameter becomes the height $H$ of the container. When $H$ is small, the heat transfer happens 
 entirely by conduction. However, as $H$ increases, such a steady state becomes unstable and 
 bifurcation of solutions occurs, see e.g. Chapter 6 in \cite{AP}. One then starts observing 
 {\em convection rolls}, where the heated fluid, becoming lighter, moves upward and then 
 returns near the bottom after cooling down driven by gravity, see some description in \cite{Chan}.  When 
 $H$ becomes larger, the 
 formation of thin conducting layers near the boundary is observed, while in the bulk of the 
 container heat transfer mostly happens via convection after a cascade of bifurcations, 
 generating chaotic dynamics and fully developed turbulence  (\cite{Kad}). 

We will consider in this paper the Boussinesq approximation in the infinite Prandtl-number limit, namely 
when the viscosity of the fluid is much bigger than the thermal diffusion.  For the derivation of the governing 
equations of the model   we still refer to \cite{Chan}. Denoting the velocity by ${\bf u} = u \, {\bf e}_1 + v \, {\bf e}_2 + w \, {\bf e}_3$,  
the temperature by $T$ and the pressure by $p$, the infinite Prandtl-number 
limit of the Boussinesq equations in a container of $\R^3$ is given by 
\begin{equation}\label{eq:Bous-1}
  \partial_t T + {\bf u} \cdot \n T - \Delta T = 0; 
\end{equation}
\begin{equation}\label{eq:Bous-2}
  \n \cdot {\bf u} = 0; 
\end{equation}
\begin{equation}\label{eq:Bous-3}
 - \Delta {\bf u} + \n p = T {\bf e}_3. 
\end{equation}
We supplement these equations  with periodicity in $(x_1,x_2) \in [0,H] \times [0,H]$, 
and Dirichlet boundary conditions, namely 
  \begin{equation}\label{eq:T-bdry}
  T = \begin{cases}
  1 & \hbox{ for } x_3 = 0; \\ 
  0 & \hbox{ for } x_3 = H, 
  \end{cases} \qquad \qquad {\bf u} = 0 \quad \hbox{ for } x_3 = 0, H. 
  \end{equation}
Using the incompressibility condition \eqref{eq:Bous-2}, one can eliminate the pressure term and obtain for the 
${\bf e}_3$-component of the velocity $w$ the equations 
\begin{equation}\label{eq:d2w}
  \Delta^2 w = - \Delta_h T; 
\end{equation}
\begin{equation}\label{eq:w-bdry}
  w = \partial_{x_3} w = 0 \qquad \qquad \hbox{ for } x_3 =  0,H, 
\end{equation}
see \cite{OS}, where $\Delta_h$ denotes the {\em horizontal Laplacian} in the 
variables $x_1, x_2$. In the sequel we will always assume $H \gg 1$.

\medspace \medspace

\noindent The following quantity measures the average vertical heat flux in terms of the  steady state conduction heat flux, 
see formula (7) in \cite{OS}. 

\begin{df}\label{d:Nu} The {\em Nusselt number} $\textup{Nu}$ is defined as 
	$$
	  \textup{Nu} = \frac{1}{H} \int_0^H \langle ({\bf u} T - \n T) \cdot {\bf e}_3 \rangle \, dx_3, 
	$$
	where 
	$$
	  \langle h \rangle(x_3) = \limsup_{t_0 \to + \infty} \frac{1}{t_0} \int_0^{t_0} 
	  \int_{I \times I} h(t,x_1,x_2,x_3) \, dx_1 dx_2 dt. 
	$$	
\end{df}
Thus the angle brackets denote the horizontal space and time averages. Recall that here we are using the container height $H$
as the scaling  parameter, following the notation in \cite{OS}.  Our main result in this paper is the following theorem. 

\begin{thm} \label{t:main}
	There exists a universal constant $C_0 > 0$ such that as $H \to + \infty$
	$$
	\textup{Nu} \leq C_0. 
	$$
\end{thm}

\medspace

\noindent In most of the physics literature 
the dimensionless parameter appears as a coefficient in \eqref{eq:Bous-3}: 
$$
  - \Delta {\bf u} + \n p = \textup{Ra} \, T {\bf e}_3, 
$$
where $\textup{Ra}$ stands for the Rayleigh number.  Formula (2.6) in \cite{DOR}  
defines the Nusselt number using this scaling, which we refer to as $\textup{Nu}_{{{\textup{phys}}}}$. 
	The theorem above corresponds in the physics language to the following result.

\begin{cor} \label{c:main}
In the limit $\textup{Ra} \to +\infty$ the following bound holds for the Nusselt number 
	$$
	  \textup{Nu}_{{\textup{phys}}} \leq C_0 (\textup{Ra})^{\frac{1}{3}}. 
	$$
\end{cor}

\medspace 

\noindent We now provide a brief history of this problem. Malkus (\cite{mal})  and Howard (\cite{How2}) gave compelling boundary-layer  arguments predicting the  scaling $\textup{Nu}_{\textup{{phys}}} \simeq  (\textup{Ra})^{\frac{1}{3}}$ (we are back to the 
notation in Definition \ref{d:Nu}). The derivation of upper bounds in convective heat transport begins with  seminal papers by Howard  (\cite{How}, \cite{How3}) and further work by Busse (\cite{Bus}, \cite{Bus2}). Some 
sub-optimal bounds, in terms of higher powers of $\textup{Ra}$, were proven in \cite{How}, \cite{DC}. 
A powerful and beautiful new idea was introduced to the subject  by Constantin and Doering, \cite{CD95}, 
and called the \emph{background field method}. To a background temperature one associates 
a quadratic form that needs to be non-negative definite: one then tries to minimize a 
suitable integral on such backgrounds (see also \cite{DC92}, \cite{DC},  \cite{DC01}). This method 
allowed the authors to obtain an optimal bound up to a logarithmic factor in $\textup{Ra}$, and  was carried forward in another  paper with a delicate analysis of singular integrals by Doering, Otto and Reznikoff, \cite{DOR}. 
More recently, in \cite{OS} another 
improvement of the estimate from the background field method was derived, showing at the same time 
some  of the limitations of such a method, see also \cite{NO}. In the same paper, a sharper 
estimate was derived using a maximum principle for the temperature, which always lies in the 
range $[0,1]$ (see also \cite{CD99}). 
More precisely,  in \cite{OS} it is shown that 
   \begin{equation}\label{eq:bs-OS-log}
   \textup{Nu} \leq C  (\textup{Ra})^{\frac{1}{3}} \log^{\frac{1}{3}} \log \textup{Ra}. 
\end{equation}
For further reading about this problem there is an extensive review article by Ahlers et al. (\cite{AGL}).

It is perhaps worthwhile to point out that numerical studies by Ierley, Kerswell and Plasting (\cite{IKP}) predicted in the 
high Rayleigh number regime, the more precise bound
$$
  \textup{Nu}_{{\textup{phys}}}  \leq 0.139 \, (\textup{Ra})^{\frac{1}{3}}. 
$$
Our methods do not give any information about the above numerical constant $0.139$. 

We also refer to 
\cite{CNO} and \cite{Wa} for related results on upper bounds in the finite-Prandtl number regime. 

\medspace \medspace

\noindent In the following we will consider the quantity
$$
  \th = T - \langle T \rangle, 
$$
which by \eqref{eq:T-bdry} satisfies the uniform bound 
\begin{equation}\label{eq:th<1}
   |\th| \leq 2\qquad \quad \hbox{ in }  \quad \{ 0 \leq x_3 \leq H \}. 
\end{equation}

We note that \eqref{eq:d2w}-\eqref{eq:w-bdry} can be rewritten as 
\begin{equation}\label{eq:d2wdht}
	\begin{cases}
		  \Delta^2 w = \Delta_h \th & \hbox{ in } \R^2 \times [0,H]; \\
		  w = \partial_{x_3} w = 0 & \hbox{ on } \{ x_3 = 0, H \}, 
	\end{cases}
\end{equation}
replacing $\theta(t,x_1,x_2,x_3)$ with $\theta(t,x_1,x_2,x_3) - \fint_{[-H/2,H/2]^2} \theta(t,x_1,x_2,x_3) dx_1 dx_2$, 
with $\fint$ denoting averaged integrals, 
we can assume that for all values of $x_3$ the function in the right-hand side of \eqref{eq:d2wdht} satisfies  
\begin{equation}\label{eq:zero-av-th}
  \int_{[-H/2,H/2]^2} \th \, dx_1 dx_2 = 0. 
\end{equation}
This does not affect the function $w$, and will be crucial for us. 

\

Taking the horizontal space- and time-average to \eqref{eq:Bous-1},  
using the boundary conditions and the constancy of vertical heat flux 
one finds that 
$$
  \textup{Nu} = \langle T \, w \rangle - \partial_{x_3} \langle T \rangle, 
$$
where $w$ is as in \eqref{eq:d2wdht}. 
In \cite{OS} (Section C, proof of Theorem 2) it was noticed  that since $\langle T \rangle \geq 0$  
and since $\sup_{x_3} \langle \th \rangle \leq \langle \sup_x \th^2 \rangle$, by the boundary conditions 
\eqref{eq:T-bdry} and \eqref{eq:w-bdry} one has 
\begin{equation}\label{eq:bd-Nu-OS}
  \textup{Nu} = \frac{1}{\d} \left( \int_0^\d \langle \th \, w \rangle \, dx_3 + 1 - \langle T|_{x_3 = \d} \rangle \right) \leq \sup_{x_3 \in (0,\d)} \langle w^2 \rangle^{\frac 12} + \frac{1}{\d} \qquad \quad 
  \hbox{ for } \d \in (0,1). 
\end{equation} 
In the same paper, the quantity $\langle w^2 \rangle$ was controlled for $\d$ small via a bound on 
$\sup_{x_3 \in (0,1)} \langle (\partial^2_{x_3} w)^2 \rangle$, leading to \eqref{eq:bs-OS-log} mostly 
via Fourier analysis, Caccioppoli-type estimates and the maximum principle.

We get instead  uniform  bounds  on $\langle \th \, w \rangle$ for $x_3 \in [0,1]$ without passing through $L^\infty$-estimates on second-order  
derivatives of $w$, but rather using Carleson measures and a variant of a theorem by Coifman and Meyer from \cite{CM}. 
Carleson measures, recalled in the Appendix, 
%
were introduced in \cite{Car}  
for the  resolution of the Corona problem. Later, in the seminal paper \cite{FS} a 
characterization of BMO functions was obtained via Carleson measures, see Theorem 
\ref{t:FS}, used in \cite{Fe}, \cite{FS} to show that the 
dual of Hardy's space is BMO.  Recall that $f \in \textup{BMO}(\R^n)$ if one has a 
uniform control of a suitable integral involving $f$ and its average for all 
balls $B \subseteq \R^n$ (see Chapter VI in \cite{St0}), and that clearly one has 
%
\begin{equation}\label{eq:in-BMO}
	  \|f\|_{\textup{BMO}(\R^n)} \leq 2 \|f\|_{L^\infty(\R^n)}. 
\end{equation}


We explain next the main steps in the 
strategy to achieve our goal. 
First, since we are dealing with a bi-Laplace equation on a large strip (with both Dirichlet and Neumann boundary 
conditions), it is useful to understand the  solution of  such a  boundary-value problem in a half-space, which represents a prototype 
situation when we are close to the bottom of the strip. 
In this case, the expression of the Green's kernel for the counterpart of \eqref{eq:d2wdht} in $\R^3_+$ is explicitly known, and we prove for it 
suitable decay properties, as well as for some if its derivatives, see Lemma \ref{l:asy-K0}. 

We then study the solution of the  problem  
\begin{equation}\label{eq:intro-lower-strip}
	\begin{cases}
		\Delta^2 \tilde{w}_1 = \chi_{\{ 0 \leq x_3 \leq H/2 \}} \Delta_h \th & \hbox{ in } \R^3_+; \\ 
		\tilde{w}_1 = 0 & \hbox{ on } \{ x_3 = 0 \}; \\  
		\partial_{x_3} \tilde{w}_1 = 0 & \hbox{ on } \{ x_3 = 0 \},
	\end{cases}
\end{equation}
that can be represented by the above kernel, forgetting about the boundary conditions 
on the upper side of the strip. Since we only care about the behavior of the solution 
for $x_3 \in [0,1]$, understanding solutions to the latter problem will be a key point in the argument. 

Using the representation formula with pointwise bounds on the kernel, that has cubic decay, would lead to 
terms diverging logarithmically in $H$, which we want to avoid. Out strategy then 
consists in estimating the leading term in the kernel, which we recognize (up to 
a suitable multiple)  as the 
horizontal Laplacian of the Poisson kernel of the upper half space, see Remark \ref{r:ker-Poi}. 
Via properties of convolutions and some integration by parts, by the above characterization of the kernel we are left with estimating 
the interaction of the leading term of $\tilde{w}_1$ with $\theta$ by means of the 
harmonic extension to  upper half-spaces of restrictions of $\theta$ to horizontal planes. 

We can then control the principal part of $\langle \th \, w \rangle$ in Proposition 
\ref{p:w1} via bilinear estimates 
from \cite{CM}, for which we present a variant in Section \ref{s:prel}, see Proposition 
\ref{p:CM-mod} for the specific statement.  Such estimates have the 
advantage to exploit the $L^2$-integrability of one of the factors, allowing us to 
localize the integration in a suitable box of size $H$.

%
%
%

%
%
%
%
%
%
%

	We then  subtract to $\tilde{w}_1$ a correction  $\hat{w}_1$ that is bi-harmonic in the 
	strip and fixes the boundary conditions on $\{ x_3 = H \}$, which were not  satisfied 
	by $\tilde{w}_1$ on $\{x_3 = H\}$: this step is  mostly worked out via Fourier analysis, more similarly 
	to \cite{OS},  needing to keep track of the uniformity of the estimates when $H$ is large. 
	This 
	gives control on a {\em surrogate} of the function $w$ in  \eqref{eq:d2wdht}, for which 
	the right hand-side is set equal to zero on the upper half-strip. 
However, the same argument works 
	when replacing in \eqref{eq:intro-lower-strip} $\chi_{\{ 0 \leq x_3 \leq H/2 \}}$ by $\chi_{\{ H/2 
		\leq x_3 \leq H \}}$: for the latter we can indeed control the $L^\infty$ norm near the lower boundary 
	due to bi-harmonicity and the local boundary conditions.  
	The estimates obtained in this way allow to obtain  
	 the conclusion of Theorem \ref{t:main}.

\

The plan of the paper is as follows. In Section \ref{s:prel} we recall some 
basic facts about Poisson's kernel and state a variant of a bilinear estimate from \cite{CM}. 
We then study the Green's function for the bi-Laplacian in the half-space, providing 
some asymptotic behaviour on its horizontal Laplacian via scaling  and invariance properties. In Section \ref{s:approx} 
we derive estimates on  solutions to the counterpart of \eqref{eq:d2wdht} in the half-space,  imposing Dirichlet and Neumann boundary conditions
on $\{x_3= 0\}$ only,  obtaining a key estimate on averaged integrals 
over horizontal periodic squares. In Section \ref{s:pf} the 
boundary conditions are then fixed on the upper layer $\{x_3 = H\}$, and  uniform bounds on Nusselt's number are derived. 
An appendix  describes the proof of Proposition \ref{p:CM-mod}, adapting the arguments for Theorem 33 in \cite{CM}.

\medspace \medspace \medspace 

\noindent  \underline{Notation}. Throughout the paper, the letter $C$ will denote a large, but fixed, 
positive constant which is allowed to vary from one formula to another, and even within the same formula. We also use the 
standard symbol $O(\cdot)$ for a quantity that is upper bounded by a fixed constant 
times its  argument.

\medspace \medspace \medspace 

\begin{center}
{\bf Acknowledgments} 
\end{center}

\noindent S.C. is grateful to Scuola Normale Superiore for the kind hospitality. 
A.M. is supported by the project {\em Geometric problems with loss of compactness} from Scuola Normale Superiore
and by the PRIN Project 2022AKNSE4 “Variational and Analytical aspects of Geometric PDE”.
He is also member of GNAMPA as part of INdAM. The authors are grateful to Christian Seis for pointing out a problem in a previous version of the paper, 
and to Camilla Nobili for useful comments.

\section{Some preliminary results}\label{s:prel}

In this section we collect some useful preliminary facts. We recall 
some properties of Poisson's kernel and a result from \cite{FS} 
concerning harmonic extensions from $\R^2$ to the half three-space, 
stating then a variant of a bilinear estimate from \cite{CM}.  
Moreover, we analyze the 
asymptotic behavior of a another kernel in $\R^3_+$, useful to represent the 
solutions of \eqref{eq:intro-lower-strip}.

\subsection{Poisson's kernel, a biharmonic kernel and a bilinear estimate}

We first recall the 
definition and some properties of Poisson's kernel in the upper half-space, which in 
three dimensions is given by the formula 
\begin{equation}\label{eq:Poisson}
	P_s(x') = \frac{1}{2 \pi} \frac{s}{\left( s^2 + |x'|^2 \right)^{\frac 32}}; \qquad \quad s > 0, \, x' \in \R^2. 
\end{equation}
If $f(x')$ is a bounded function on $\R^2$, then 
\begin{equation}\label{eq:Fs}
	  F(x',s) := (P_s * f)(x) 
\end{equation}
is the harmonic extension of $f$ to $\R^3_+$. For later convenience, 
we notice that 
\begin{equation}\label{eq:Lapl-Poisson}
	\Delta_{h} P_s(x') = \frac{3 s \left( 3 |x'|^2 - 2 s^2\right)}{2 \pi  \left(s^2+ |x'|^2\right)^{\frac 72}}.
\end{equation}
We also recall the following result (in a particular case),  giving a characterization of BMO functions.

\begin{thm}[\cite{FS}, Theorem 3]  \label{t:FS}
	There exists a constant $C_0 > 0$ such that, for 
	 $f(\cdot) \in \textup{BMO}(\R^2)$ and letting 
	$$
	  T(x^0,h) := \left\{ (x,s) \in \R^3_+ \; : \; |x - x^0| \leq h, s \in (0,h) \right\}; \qquad 
	  x^0 \in \R^2, \, h > 0, 
	$$
	one has  
	\begin{equation} \label{eq:int-cube}
		\sup_{x^0 \in \R^2} \int_{T(x^0,h)}  |t \nabla_h F(x',t)|^2 \, d x' \, \frac{dt}{t} \leq C_0  \, h^2 \|f\|_{\textup{BMO}(\R^2)}^2,
	\end{equation}
where $F$ is as in \eqref{eq:Fs} amd $\nabla_h$ stands for the horizontal gradient. 
\end{thm}

\begin{rem}\label{r:FS}
(i) Even though the quadratic dependence in the \textup{BMO} norm  is not explicitly stated in 
Theorem 3 of \cite{FS}, it appears there in the second formula on page 147. 

(ii)	Provided  $(P_1 * f)(0) \in \R$,   finiteness of the integral in \eqref{eq:int-cube} implies in turn that 
$f$ is of class \textup{BMO}. 
\end{rem}

\

\begin{figure}[!htb]
\begin{center}
	\minipage{0.60\textwidth}
	\includegraphics[width=\linewidth]{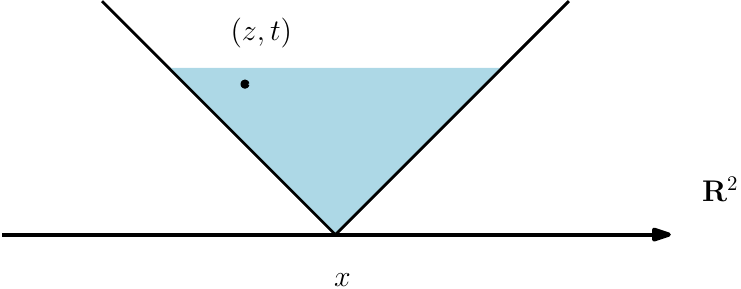}
	\caption{\small the non-tangential maximal function}\label{fig:ntm}
	\endminipage\hfill
\end{center} 
\end{figure}

For a measurable and bounded function $F(z,t)$ defined on $\R^3_+$, the 
{\em non-tangential maximal function} (see Figure \ref{fig:ntm}) $N F : \R^2 \to \R$ is given by 
\begin{equation}\label{eq:NTM}
  N F(x) = \sup \left\{ |F(z,t)| \; : \;  |x-z| < t \right\}. 	
\end{equation}
For $F$ as before and $t > 0$, we define $F^t : \R^2 \to \R$ as 
$$
  F^t(x') := F(x',t). 
$$

The above function is useful in the proof 
of \cite[Theorem 33]{CM}, 
for which we have the following variant.

\begin{pro}\label{p:CM-mod}
%
	Let $\psi$ be the vectorial  convolution kernel in $\R^2$ given by 
	$$
	\psi(x') = \nabla_{x'} P_1(x'). 
	$$ 
	For $t > 0$ let $\psi_t(\cdot) = t^{-2} \psi(t^{-1} \cdot)$. Let $\rho$ be a 
\textup{BMO} function in $\R^2$ and let $F$ be a bounded function on $\R^3_+$ 
with compact support. Let also $m : (0,+\infty) \to \R$ be measurable with 
$|m(t)| \leq 1$ for all $t$. Then there exists a  
constant $C > 0$, independent of $\rho, F$ and $m$, such that letting $\check{F}(x',t) = (\psi_t* F^t) (x')$, one has  
$$
\left( \int_{\R^2} \left| \int_0^\infty  (\psi_t * \rho)(x') \cdot (\psi_t* F^t) (x')  m(t) \frac{dt}{t}  \right|^2  dx'  \right)^{\frac 12}
\leq C \| \rho \|_{\textup{BMO}} \|N \check{F} \|_{L^2(\R^2)}.
$$		
\end{pro}

\begin{rem}\label{r:CM-gen}
	The same conclusion holds for a convolution kernel verifying 
	$$
	\int_{\R^2} \psi (x') dx' 	= 0, 
	$$
	and such that its Fourier transform satisfies  
	\begin{equation}
		\begin{cases} \label{eq:111-0}
			| \partial^\alpha_\xi \hat{\psi}(\xi)| \leq C(\alpha) e^{- C(\alpha)^{-1} |\xi|} & \hbox{ for } |\xi| \geq 1; \\ 
			|\xi|^\alpha |\partial^\alpha_\xi  \hat{\psi}(\xi)| \leq C(\alpha) |\xi| & \hbox{ for } |\xi| \leq  1. 
		\end{cases}
	\end{equation}	
In fact, via the same proof, one can also replace $\psi_t* F^t$ in Proposition \ref{p:CM-mod} with $\varphi_t* F^t$ for a kernel $\varphi$ 
as in \eqref{eq:111-0} but without vanishing integral.

We will prove  in the appendix that $\|N \check{F} \|_{L^2(\R^2)}$  is finite for the  
kernel in Proposition \ref{p:CM-mod}. This property can be indeed verified 
under the previous more general assumptions. 
\end{rem}

The difference of Proposition \ref{p:CM-mod} compared to the original version is that we allow a dependence in the 
$t$-variable in the second convolved function. The argument follows the original one, 
and in the appendix we describe how the proof in \cite{CM} adapts to this 
situation as well.

The result applies in particular to the (vectorial) case of the horizontal gradient of Poisson's kernel 
$$
	 \psi(x') = \nabla_h P_1(x'), 
	$$
where $P_1$ is as in \eqref{eq:Poisson} with $s  = 1$.

%
%
%
%
%

\subsection{A useful kernel in the upper half-space}

We are interested next in deriving some asymptotic estimates involving the Green's function for the bi-Laplacian with zero Dirichlet and Neumann boundary conditions in $\R^3_+ = \left\{ (\eta_1,\eta_2,\eta_3) \; : \; \eta_3 > 0 \right\}$, 
namely the solution $G(\xi,\eta)$ (decaying to zero as $\eta_3 \to + \infty$) of 
$$
  \begin{cases}
  \Delta^2_{\eta} G(\xi,\cdot) = \delta_{\xi} & \hbox{ in } \R^3_+; \\ 
  G(\xi,\cdot) = 0  & \hbox{ on } \partial \R^3_+; \\ 
  \partial_{\eta_3} G(\xi,\cdot) = 0 & \hbox{ on } \partial \R^3_+. 
  \end{cases}
$$
Recall that from Lemma 8 in \cite{GR}, $G$ has the exact expression 
$$
  G(\xi,\eta) = \frac{1}{16 \pi} |\xi-\eta| \int_{1}^{\frac{|\xi^*-\eta|}{|\xi-\eta|}} (t^2-1) \, t^{-2} dt, 
$$
where $\xi^*$ stands for the reflection of $\xi$ across the boundary, $\{\eta_3 = 0 \}$. 
From an explicit computation one  obtains the following formula for $G$ 
$$
  G(\xi,\eta) =  \frac{1}{16 \pi} |\xi-\eta| \left[ \frac{|\xi^*-\eta|}{|\xi-\eta|} + \frac{|\xi-\eta|}{|\xi^*-\eta|} -2 \right], 
$$
which after further manipulation becomes 
$$
  G(\xi,\eta) = \frac{1}{16 \pi} \frac{\left( |\xi^*-\eta| - |\xi-\eta| \right)^2}{|\xi^*-\eta|}. 
$$
Computing the horizontal Laplacian $\Delta_h^\eta = \frac{\pa^2}{\pa \eta_1^2} + 
\frac{\pa^2}{\pa \eta_2^2}$, one then finds 
\begin{equation}\label{eq:formula-Green}
 K(\xi,\eta) := \Delta_h^\eta G(\xi,\eta) =  \frac{1}{8 \pi} \left[ \frac{1}{|\xi^*-\eta|} - \frac{1}{|\xi-\eta|} 
-  \frac{(\xi_3-\eta_3)^2}{|\xi-\eta|^3} + \frac{\xi_3^2+\eta_3^2}{|\xi^*-\eta|^3} 
+ 6 \frac{\xi_3 \eta_3 (\xi_3+\eta_3)^2}{|\xi^*-\eta|^5}  \right]. 
\end{equation}
%
We observe  the following three 
covariance properties for all $\xi, \eta \in \R^3_+$ and for all $v' \in \R^2$: 
\begin{equation}\label{eq:symm-K}
  K(\l \xi, \l \eta) = \frac{1}{\l} K(\xi,\eta), \quad \l > 0; \qquad \quad K(\xi,\eta) = K(\eta,\xi); \qquad \quad
  K(\xi + (v',0), \eta + (v',0)) = K(\xi,\eta). 
\end{equation}
By these, it will be sufficient to understand the asymptotics for $|\eta| \to + \infty$ of 
\begin{equation}\label{eq:K0}
   K_0(\eta) = K(\xi_0,\eta); \qquad   \qquad \xi_0 = (0,0,1). 
\end{equation}
We also notice that  the function $K_0$ is even with respect to  reflection around the $\eta_3$-axis, namely   
\begin{equation}\label{eq:K0-refl}
K_0(\eta_1, \eta_2, \eta_3) = K(- \eta_1, - \eta_2, \eta_3). 
\end{equation}
%
%

We will now deduce some asymptotics on the function $K_0$. 

\begin{lem}\label{l:asy-K0}
	For $\xi_0$ and $K_0$ as in \eqref{eq:K0}, we have the upper bound 
	\begin{equation}\label{eq:est-short-dist-0}
		|K_0(\eta)| \leq \frac{C}{|\eta-\xi_0|}; \qquad \quad |\eta-\xi_0| \leq 8. 
	\end{equation}	
	Moreover, for a vector $(a,b,c) \in \R^3$ of unit norm, $c \geq 0$, and for $R > 0$ large one has 
	\begin{equation}\label{eq:K0K0}
			K_0(Ra, Rb,  Rc) = \frac{3 c^2 \left(3 a^2 +3 b^2 -2 c^2\right)}{4 \pi  R^3}+O\left( R^{-4} \right), 
	\end{equation}
		\begin{equation}\label{eq:dec-1st-der}
		\nabla_{\eta} K_0(Ra, Rb,  Rc) = O(R^{-4});
	\end{equation}
		\begin{equation}\label{eq:dec-2nd-der}
		\nabla^{(2)}_{\eta,\eta} K_0(Ra, Rb,  Rc) = O(R^{-5}). 
	\end{equation}
	The above estimates are uniform as $R \to + \infty$ in the choice of the unit vector $(a,b,c)$, with $c \geq 0$. Furthermore, 
	there exists $C > 0$ such that
	\begin{equation}\label{eq:basso}
			K_0(\eta',\eta_3) \leq \frac{C}{|\eta'|^3}; \qquad \eta' \in \R^2, |\eta'| \geq 1, \quad \eta_3 \in (0,1]. 
	\end{equation}
\end{lem}

\begin{pf}
The proof of \eqref{eq:K0K0} can be deduced from some cancellations after Taylor expansions of the following 
quantities, corresponding to each of  term in \eqref{eq:formula-Green}, choosing $\eta = R(a,b,c)$, with $a^2 + b^2 + c^2 = 1$, 
$\xi = (0,0,1)$ and $\xi^* = (0,0,-1)$: 
$$
  \frac{1}{(R^2 a^2+R^2 b^2+(R c+1)^2)^{\frac 12}} = \frac{1}{R } -\frac{c}{R^2 } + \frac{2 c^2-a^2-b^2}{2 R^3 }+O\left(\frac{1}{R^4}\right); 
$$
$$
 \frac{1}{\sqrt{R^2 a^2+R^2 b^2+(R c-1)^2}} = \frac{1}{R }+\frac{c}{R^2 }+\frac{2 c^2-a^2-b^2}{2 R^3 }+O\left(\frac{1}{R^4}\right); 
$$
$$
\frac{(R c-1)^2}{\left(R^2 a^2+R^2 b^2+(1-R c)^2\right)^{3/2}} = 
\frac{c^2}{R }+\frac{c^3-2 c \left(a^2+b^2\right)}{R^2 }+\frac{2 c^4-11 c^2 \left(a^2+b^2\right)+2 \left(a^2+b^2\right)^2}{2 R^3 }+O\left(\frac{1}{R^4}\right); 
$$
$$
\frac{R^2 c^2+1}{\left(R^2 a^2+R^2 b^2+(R c+1)^2\right)^{3/2}} = 
\frac{c^2}{R }-\frac{3 c^3}{R^2 }+\frac{c^2 \left(a^2+b^2\right)+2 \left(a^2+b^2\right)^2+14 c^4}{2 R^3 }+O\left(\frac{1}{R^4}\right); 
$$
$$
\frac{R c (R c+1)^2}{\left(R^2 a^2+R^2 b^2+(R c+1)^2\right)^{5/2}}
= \frac{c^3}{R^2 } +  \frac{2 c^2 \left(a^2+b^2\right)-3 c^4}{R^3 }+O\left(\frac{1}{R^4}\right). 
$$
The proofs of the gradient and the hessian estimates can be obtained in a similar manner.

While the condition $c \geq 0$ is needed in \eqref{eq:K0K0} and the next two formulas to have the argument
of $K_0$ in $\R^3_+$, the four estimates hold regardless of the    
sign of $c$, allowing to show \eqref{eq:basso} as well. 
\end{pf}

We observe that, by a Talor expansion, an identical estimate as in \eqref{eq:K0K0} and 
in the next two formulas holds for $K_0(\xi_0 + R(a,b,c))$ as well.

\

Using the above  result, we obtain next the following properties of $K(\cdot, \cdot)$. 

%

\begin{pro}
There exists $C > 0$ such that for $x_3, y_3  > 0$, $x \neq y$, the  kernel $K$ satisfies 
\begin{equation}\label{eq:est-short-dist}
  |K(x,y)| \leq \frac{C}{|x-y|}, \quad \qquad \hbox{for } |x - y| \leq 8 x_3. 
\end{equation}	
Moreover, for any  unit vector  $(a,b,c)$ with $c \geq 0$, we have,  
uniformly in $c \geq 0$ for $R > 2 x_3$: 
\begin{equation} \label{l:2-1}
 K(x, x+ R(a,b,c))  =   
\frac{3 x_3^2 c^2  \left(3 a^2+3 b^2-2 c^2\right)}{4 \pi R^3 }
+O\left(\frac{x_3^3}{R^4}\right); 
\end{equation}
\begin{equation}\label{l:2-3}
	(\nabla_{\eta} K) (x,x+R(a,b,c))  =   
	O\left(\frac{x_3^2}{R^4}\right), 
\end{equation}
and   
\begin{equation}\label{l:2-4}
		(\nabla^{(2)}_{\eta, \eta} K) (x,x+R(a,b,c))  =   
	O\left(\frac{x_3^2}{R^5}\right). 
\end{equation}
\end{pro}

\begin{pf}
	We only show \eqref{l:2-1}, the other two estimates begin similar. 
	
	From the third and first properties in \eqref{eq:symm-K} we find 
	$$
	 K(x, x+ R(a,b,c))  =   K((0,x_3),(0,x_3)+ R(a,b,c)) = \frac{1}{x_3} K((0,0,1),(0,0,1)+R/x_3(a,b,c)). 
	$$
	Recalling \eqref{eq:K0}, this becomes 
	$$
	 K(x, x+ R(a,b,c))  =  \frac{1}{x_3} K_0(\xi_0 + R/x_3(a,b,c)). 
	$$
	From \eqref{eq:K0K0} and the observation after Lemma \ref{l:asy-K0} we then obtain 
	$$
	K(x, x+ R(a,b,c))  =  \frac{1}{x_3} \left(  \frac{3 c^2 \left(3 a^2 +3 b^2 -2 c^2\right)}{4 \pi}  
	\frac{x_3^3}{R^3} +O\left( \frac{x_3^4}{R^4} \right)   \right),  
	$$
	as desired. 
\end{pf}

\begin{rem}\label{r:ker-Poi}
	We notice  that, comparing with \eqref{eq:Lapl-Poisson}, for $y_3 > 2 x_3$  
	$$
	 K(x,y) = \frac{1}{2} x_3^2 y_3 (\Delta_h P_{y_3}) (x'-y') + O\left(\frac{x_3^3}{|x-y|^4}\right).
	$$
In fact, from a Taylor expansion of \eqref{eq:Lapl-Poisson} and \eqref{l:2-1} one finds that 
$$
 y_3  (\Delta_h P_{y_3}) (x'-y') =  \frac{3 y_3^2 (3 |x'-y'|^2 - 2y_3^2)}{2 \pi(y_3^2 + |x'-y'|^2)^{\frac{7}{2}}} 
  =  \frac{3 (y_3-x_3)^2 (3 |x'-y'|^2 - 2(y_3-x_3)^2)}{2 \pi((y_3-x_3)^2 + |x'-y'|^2)^{\frac{7}{2}}} 
  +  O \left( \frac{ x_3 }{|x-y|^4} \right).
$$	
	This observation will be crucial for the estimates in the next section. 
\end{rem}

\section{Estimates on approximate solutions to \eqref{eq:d2wdht} }\label{s:approx}


We now consider two cut-off functions in the vertical variable  $\chi_{\{ 0 \leq x_3 \leq H/2 \}}, \chi_{\{ H/2 \leq x_3 \leq H \}}$ and split 
$\theta$ as $\theta = \theta_1 + \theta_2$, with 
\begin{equation}\label{eq:th-i}
\theta_1 = \theta \chi_{\{ 0 \leq x_3 \leq H/2 \}}; \qquad \quad 
\theta_2 = \theta \chi_{\{ H/2 \leq x_3 \leq H \}}. 
\end{equation}
For $K$ as in \eqref{eq:formula-Green}, we will prove in this section a key estimate involving the function 
\begin{equation}\label{eq:repr}
  \tilde{w}_1(x) := \int_{\R^3_+} K(x,y) \, \th_1(y) dy 
\end{equation}
near the bottom  of the strip $\{0 \leq x_3 \leq H\}$. We begin with a preliminary lemma, 
for which we recall Remark \ref{r:ker-Poi}.

\begin{lem}\label{l:I+II}
	There exists a fixed constant $C > 0$ such that, writing 
	  $\tilde{w}_1(x)  = I + II$ with 
	  $$
	    I = \int_{\{0 < y_3 \leq 4 x_3\}} K(x,y) \, \th_1(y) dy; \qquad 
	    \quad  II = \int_{\{4 x_3 < y_3 \leq H/2\}} K(x,y) \, \th_1(y) dy,
	  $$
	  one has the estimates 
\begin{equation} \label{eq:I-II}
	|I| \leq C x_3^2; \quad \left| II  - \frac{x_3^2}{2}  \int_{\{4 x_3 < y_3 \leq H/2\}} y_3 ( (\Delta_{h} P_{y_3}) * \th_1^{y_3} )(x')   dy_3 \right| 
	\leq  C x_3^2; \qquad  \quad x_3 \in (0,1]. 
\end{equation}
In the latter formula, the convolution is taken with respect to the horizontal variables, and $\th_1^{x_3} : \R^2 \to \R$ 
denotes the restriction of $\theta_1$ to the horizontal plane at height $x_3$, namely  
$$
 \th_1^{x_3} (x') := \theta_1(x',x_3). 
$$
\end{lem}

\begin{pf}
	To prove the first inequality in \eqref{eq:I-II}, we consider the cube $Q_x$ centered at $(x',2x_3)$ and of size $4 x_3$, 
	so its lower face lies on $\{x_3 = 0\}$, see Figure \ref{fig:cubes}.  
	\begin{figure}[!htb]
		\begin{center}
			\minipage{0.70\textwidth}
			\includegraphics[width=\linewidth]{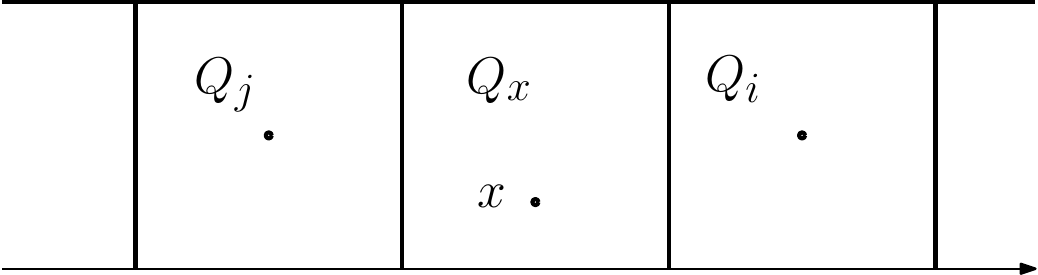}
			\caption{\small the cubes $Q_x$ and $(Q_i)_i$}\label{fig:cubes}
			\endminipage\hfill
		\end{center} 
	\end{figure}
	We next partition the strip $\{ 0 < y_3 < 4 x_3 \}$ into $Q_x$ and a 
	sequence of cubes defined as follows. For $z_i \in \Z^2 \setminus \{(0,0)\}$, 
	let $Q_i$  denote the cube of size $4 x_3$ with center at the point $p_i := (x_1,x_2, 4 x_3 z_i)$. 
	
	We then notice that, by \eqref{eq:est-short-dist} and \eqref{eq:symm-K} 
	$$
	  \left| \int_{Q_x} K(x,y) \theta_1(y) dy \right| \leq C \left| \int_{Q_x} \frac{1}{|x-y|} dy \right| 
	  \leq C \int_0^{8 x_3} \frac{r^2}{r} dr \leq C x_3^2.  
	$$
	On the other hand, concerning $Q_i$, since the decay of $K$ from \eqref{l:2-1}  is cubic towards infinity,   we have that 
	$$
	  \left| \int_{Q_i} K(x,y) \theta_1(y) dy \right| \leq C |Q_i| \sup_{y \in Q_i} |K(x,y)| 
	  \leq C x_3^3  \frac{x_3^2}{|x-p_i|^3}. 
	$$
	Therefore, by the choice of $p_i$ and $Q_i$ we have that 
	$$
	  \left|  \int_{\{ 0 < y_3 < 4 x_3 \} \setminus Q_x} K(x,y) \theta_1(y) dy \right|  \leq C 
	  x_3^5 \sum_i \frac{1}{|x-p_i|^3} \leq C x_3^5 \int_{x_3}^{\infty} \frac{r}{r^3} dr \leq C x_3^2.  
	$$
	The last three formulas prove that $|I| \leq C x_3^2$, as desired. 
	
	\ 
	
	The estimate of $II$ simply follows from Remark \ref{r:ker-Poi}, that yields 
	\begin{eqnarray*}
		\left| II  - \frac{x_3^2}{2}  \int_{\{4 x_3 < y_3 \leq H/2\}} y_3 ( (\Delta_{h} P_{y_3}) * \th_1^{y_3} )(x')   dy_3 \right|  
		& \leq & C \int_{\{4 x_3 < y_3 \leq H/2\}} \frac{x_3^3}{|x-y|^4} dy \\ 
		& \leq & C x_3^3 \int_{3 x_3}^\infty \frac{s^2}{s^4} ds \leq C x_3^2, 
	\end{eqnarray*}
concluding the proof. 
\end{pf}

We next estimate the principal part of $II$ appearing in  \eqref{eq:I-II}, integrated versus 
$\theta$ in $x' \in [-H/2,H/2]$ and in $x_3 \in [0,1]$. 
In the next formula and in the proof below, we will sometimes write the differentials right after the corresponding 
integral signs: even though this is unconventional, we hope it will facilitate tracking the 
various domains of integration.

\begin{pro}\label{p:FS}
	Let $\tilde{w}_1$ be as in \eqref{eq:repr}. Then there exists a constant $C_1 > 0$
	\underline{independent of $H$} (large) such that for any $x_3 \in [0,1]$ 
	one has the upper bound 
\begin{equation*}
	\left|  \frac{1}{H^2} \int_{[-H/2,H/2]^2} dx' \theta_1^{x_3}(x') \, x_3^2 
	 \int_{\{4 x_3 < y_3 \leq H/2\}}  y_3 \, dy_3  ( (\Delta_{h} P_{y_3}) * \th_1^{y_3} )(x')   \right| \leq C_1. 
\end{equation*}
\end{pro}

\begin{pf}
	We first localize $\theta_1^{y_3} : \R^2 \to \R$ onto the box $[-2H,2H]^2$ and its complement in the 
	$(x_1,x_2)$-coordinates. We write 
	$$
	\theta_1^{y_3} = g_1^{y_3} + g_2^{y_3}, \qquad \quad \hbox{ with} \quad 
	g_1^{y_3}  = \theta_1^{y_3} \chi_{[-2H,2H]}, \quad 
	g_2^{y_3}  = \theta_1^{y_3} \chi_{\R^2 \setminus [-2H,2H]}. 
	$$
	We note next that for $x' \in [-H/2,H/2]$, by \eqref{eq:Lapl-Poisson} and the fact that $\|\theta\|_{L^\infty} \leq 1$, one has 
	$$
	|y_3 ((\Delta_h P_{y_3})* g_2^{y_3})(x')| \leq C \int_{\R^2 \setminus [-2H,2H]^2} 
	\frac{|\theta_1^{y_3}(w)| dw}{(y_3^2 + |x'-w|^2)^{\frac{3}{2}}} \leq C \int_{\{|w| \geq H\}} \frac{dw}{|w|^3} \leq \frac{C}{H}. 
	$$
	This implies that 
	\begin{equation*}
		\left| \int_{[-H/2,H/2]^2} \theta_1^{x_3}(x') \int_{\{4 x_3 \leq y_3 \leq H/2\}} 
		x_3^2 y_3 ((\Delta_h P_{y_3}) * g_2^{y_3})(x') dy_3 dx' \right| \leq 
		\int_{[-H/2,H/2]^2}  \int_0^H dy_3 \frac{C}{H} dx' \leq C H^2. 
	\end{equation*}

	Therefore, we will now focus on $g_1^{y_3}$. By the semi-group property of the 
	Poisson kernel we have that $P_{y_3} = P_{y_3/2} * P_{y_3/2}$, therefore using the symmetry 
	and  	associativity of   
	convolutions together with  the fact that for $f \in L^\infty(\R^2)$ 
	\begin{equation*}
		\frac{\partial}{\partial x'_i} (P_{s} * f) = \left(  \frac{\partial}{\partial x'_i} P_{s} \right) * f; 
		\qquad i = 1, 2, s > 0, 
	\end{equation*}
	we find 
	\begin{eqnarray*}
		\int_{[-H/2,H/2]^2} dx' \theta_1^{x_3}(x') \int_{\{4 x_3 \leq y_3 \leq H/2\}} 
		x_3^2 y_3 ((\Delta_h P_{y_3}) * g_1^{y_3})(x') dy_3 dx'  \\ 
		= \int_{\R^2} \int_{\{4 x_3 \leq y_3 \leq H/2\}}  (P_{y_3/2} * \left( \theta_1^{x_3} \chi_{[-H/2,H/2]^2} \right))(x') 
		x_3^2 y_3 \Delta_h (P_{y_3/2} * g_1^{y_3})(x') dy_3. 
	\end{eqnarray*}
	Integrating by parts in the horizontal variables and setting 
	$$
	\tilde{\theta}_1^{x_3} = \theta_1^{x_3}(x') \chi_{[-H/2,H/2]^2}(x'), 
	$$
	we find that  the latter integral expression is equal to  $A+B$, 
	where  
	$$
	A = \int_{\R^2 \setminus [-8H,8H]^2}  \int_{\{4 x_3 \leq y_3 \leq H/2\}} x_3^2 y_3 (  (\nabla_h P_{y_3/2}) * \tilde{\theta}_1^{x_3} )(x') \cdot 
	 ((\nabla_h P_{y_3/2}) * g_1^{y_3})(x') dy_3 dx'; 
	$$
	$$
	B = \int_{ [-8H,8H]^2}  \int_{\{4 x_3 \leq y_3 \leq H/2\}} x_3^2 y_3 ( (\nabla_h P_{y_3/2}) *   \tilde{\theta}_1^{x_3} )(x') 
	\cdot  ((\nabla_h P_{y_3/2}) * g_1^{y_3})(x') dy_3 dx'. 
	$$
	We claim that $A$ is of order $O(H^2)$: in fact, we notice that from the expression of $P_s$ in \eqref{eq:Poisson} and 
	an elementary inequality  
	$$
	|\nabla_h P_{y_3/2}(w)| \leq \frac{C}{(y_3^2 + |w|^2)^{\frac{3}{2}}}; \qquad 
	|y_3 \nabla_h P_{y_3/2}(w)| \leq \frac{C y_3}{(y_3^2 + |w|^2)^{\frac{3}{2}}}. 
	$$
	Hence, for $x' \in \R^2 \setminus [-8H,8H]^2$ we have that 
	$$
	\left|(  (\nabla_h P_{y_3/2}) * \tilde{\theta}_1^{x_3} )(x') \right| \leq 
	\int_{\R^2} \frac{|\tilde{\theta}_1^{x_3}| dw}{(y_3^2 + |x'-w|^2)^{\frac{3}{2}}} 
	\leq C \frac{H^2}{|x'|^3}. 
	$$ 
	On the other hand, still for 	 $x' \in \R^2 \setminus [-8H,8H]^2$ we obtain  
	$$
	\left|y_3 ((\nabla_h P_{y_3/2}) *  g_1^{y_3}) (x') \right| \leq 
	\int_{\R^2} \frac{|g_1^{y_3}| y_3 dw}{(y_3^2 + |x'-w|^2)^{\frac{3}{2}}} 
	\leq C \frac{H^2}{|x'|^2}. 
	$$
	Hence we find that 
	\begin{equation}\label{eq:A}
		|A| \leq C \int_{\{|x'| \geq 8 H\}} x_3^2 \int_0^{H/2} dy_3 \frac{H^4}{|x'|^5} dx' \leq C H^2. 
	\end{equation}

	We next estimate the term $B$: using the Cauchy-Schwartz inequality in $x'$ we get that 
	\begin{eqnarray*} \nonumber
		|B| & \leq & x_3^2 \left( \int_{ [-8H,8H]^2}  dx'\right)^{\frac{1}{2}} \\ 
		& \times & \left[ \int_{ [-8H,8H]^2}  \left| 
		\int_0^\infty y_3 ((\nabla_h P_{y_3/2}) * \tilde{\theta}^{x_3}_1)(x') \cdot 
		y_3 ((\nabla_h P_{y_3/2}) * g_1^{y_3})(x') \chi_S(y_3) \frac{dy_3}{y_3} \right|^2 dx' \right]^{\frac{1}{2}}, 
	\end{eqnarray*} 
	where we defined the set $S$ as 
	$$
	S = \left\{ y_3 \; | \; 4 x_3 \leq y_3 \leq \frac{H}{2} \right\}.
	$$
	Setting now 
	$$
	t = y_3; \qquad \psi_t = t \nabla_h P_{t/2}; \qquad m(t) = \chi_S(t), 
	$$
	we can write that 
	$$
	|B| \leq C H \left[ \int_{ \R^2}  \left| 
	\int_0^\infty (\psi_t * \tilde{\theta}^{x_3}_1)(x') \cdot 
	(\psi_t * g_1^{t})(x') \, m(t) \frac{dt}{t} \right|^2 dx' \right]^{\frac{1}{2}}	 
	x_3^2. 
	$$
	Noting that (in vectorial sense) $\int_{ \R^2} \psi_t(x') dx' = 0$, 
	we have by  Theorem \ref{t:FS} that 
	$$
	|(\psi_t * \tilde{\theta}^{x_3}_1)(x')|^2 \frac{dx' dt}{t}
	$$
	is a Carleson measure, see Definition \ref{df:carleson} and Theorem \ref{t:FS}, since $\tilde{\theta}^{x_3}_1(\cdot)$ is in 
	$L^\infty$ \underline{with $x_3$ fixed}. 
	Moreover, we have clearly 
	$$
	\|m(t)\|_\infty 
	\leq 1. 
	$$
	Recalling the notion of non-tangential maximal function from \eqref{eq:NTM} 
	and choosing 
	$$
	F(z,t) =g_1^t(z); \qquad  \check{F}(z,t) = (\psi_t*g_1^t)(z),
	$$ 
	applying Proposition \ref{p:CM-mod}	and \eqref{eq:in-BMO} 
	we find that
	\begin{equation}\label{eq:B}
		|B| \leq C  H \|N \check{F}\|_{L^2(\R^2)} \|\theta_1\|_{L^\infty}. 
	\end{equation}
	To estimate the above $L^2$-norm we first notice that 
	\begin{equation}\label{eq:11}
		|\check{F}(z,t)| = |(\psi_t* g_1^t)(z)| \leq C \int_{\R^2} \frac{|g_1^t(w)| t}{(t^2 + |z-w|^2)^{\frac 32}} dw \leq 
		C \int_{ \R^2} \frac{t}{(t^2 + |w|^2)^{\frac 32}} dw \leq C. 
	\end{equation}
	We now specialize to the case in which $x' \in \R^2 \setminus [-10H, 10H]^2$: 
	notice that for $|z-x'| \leq t$ we have 
	$$
	t^2 + |x'-w|^2 \leq 2 (t^2 + |z-w|^2), 
	$$
	which implies 
	$$
	P_t(z-w) \leq C P_t(x' - w); \qquad \hbox{ for } \quad |z-x'| \leq t. 
	$$
	Using this inequality we get 
	$$
	|\check{F}(z,t)| \leq C \int_{\R^2} \frac{|g_1^t(w)| t}{(t^2 + |x'-w|^2)^{\frac 32}} dw; 
	\qquad |z-x'| \leq t, \; x' \in \R^2 \setminus [-10H, 10H]^2. 
	$$
	Since we only integrate for $w \in [-2H, 2H]^2$ in the support of $g_1^t$, we then obtain 
	$$
	|\check{F}(z,t)| \leq \frac{C}{|x|^2} \int_{\R^2} |g_1^t(w)| dw \leq C \frac{H^2}{|x|^2},  
	\qquad |z-x'| \leq t, \; x' \in \R^2 \setminus [-10H, 10H]^2. 
	$$
	Recalling \eqref{eq:NTM}, combining \eqref{eq:11} and the latter estimate we find that 
	$$
	\|N \check{F}\|_{L^2(\R^2)} \leq C H. 
	$$
	From the last inequality, \eqref{eq:A} and \eqref{eq:B} we then get the desired conclusion. 
\end{pf}

From Lemma \ref{l:I+II} and Proposition \ref{p:FS} we derive the following consequence.

\begin{pro} \label{p:w1}
	Let $\th_1$ be as in \eqref{eq:th-i} and $\tilde{w}_1$ as in \eqref{eq:repr}. 
	Then there exists a fixed constant $C > 0$ independent of $H$ (large) 
	such that 
	$$
	  \left|  \frac{1}{H^2} \int_{[-H/2,H/2]^2}  \theta_1(x',x_3) \tilde{w}_1(x',x_3) dx'    \right| \leq C; 
	  \qquad x_3 \in (0,1]. 
	$$
\end{pro}


	We next consider the above function $\tilde{w}_1$ for $x_3$ close to $H$, where we have better 
	regularity properties due to its bi-harmonicity.

\begin{pro}\label{p:upper-bdry}
	There exist fixed positive constants $H_0$ and $C_0$ such that for $H \geq H_0$ one has 
	\begin{equation}\label{eq:estx3L}
	|\tilde{w}_1(x)| \leq   C_0 H^2;  \quad |\nabla \tilde{w}_1(x)| \leq C_0 \, H; \quad |\nabla^{(2)} \tilde{w}_1(x)| \leq C_0, 
	\qquad  x = (x',x_3), \; \;  x_3 \in [3/4 \, H, 5/4 \, H].
	\end{equation}
\end{pro}

\begin{pf} 
We proceed similarly to the proof of the first estimate in \eqref{eq:I-II}. 
 Assuming w.l.o.g. that $x' = 0$, we let $\tilde{Q}_i$  denote the cube of size $H/2$ centered at 
$\hat{z}_i := (H/2 \; \tilde{z}_i, H/4)$, with $\tilde{z}_i \in \Z^2$, 
see Figure \ref{fig:cubes-2}.

\begin{figure}[!htb]
	\begin{center}
		\minipage{0.75\textwidth}
		\includegraphics[width=\linewidth]{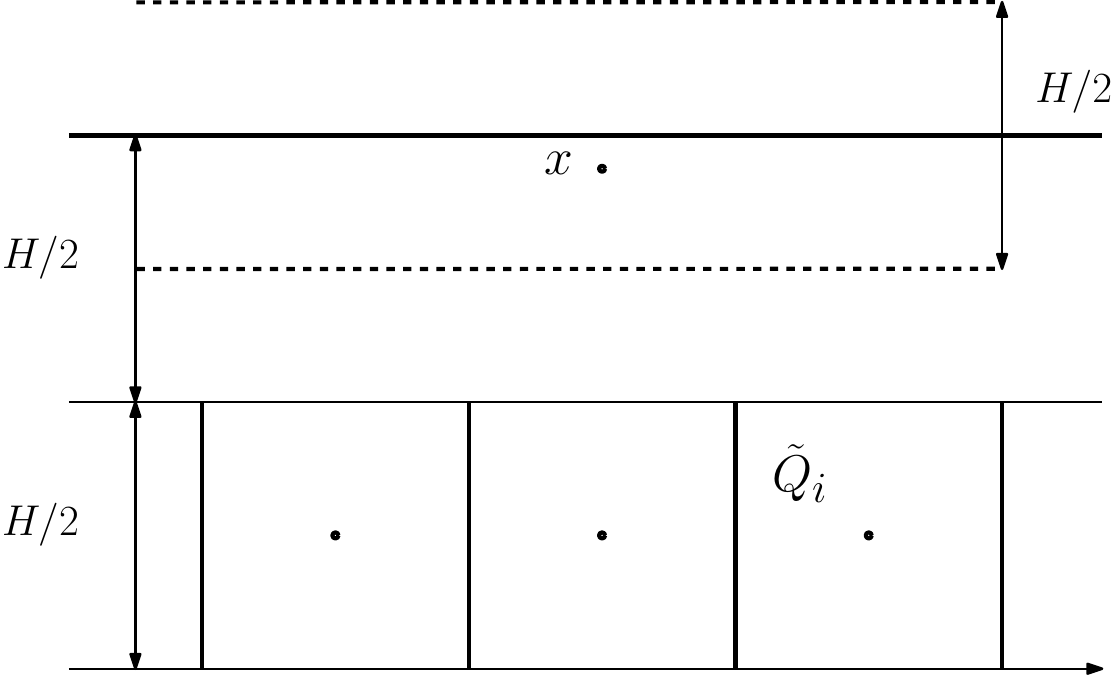}
		\caption{\small the cubes $(\tilde{Q}_i)_i$}\label{fig:cubes-2}
		\endminipage\hfill
	\end{center} 
\end{figure}

Using  \eqref{l:2-1} and the second property in \eqref{eq:symm-K}  we have that, for $x' = 0$ and $x_3 \in [3/4 \, H, 5/4 \, H]$, 
\begin{equation}\label{eq:c14}
	 \left| \int_{\tilde{Q}_i} K(x,y) \tilde{\theta}_1(y) dy  \right| =  
	\left| \int_{\tilde{Q}_i} K(y,x) \tilde{\theta}_1(y) dy  \right| \leq C \left|\int_{\tilde{Q}_i} y_3^2 dy\right|  
	\frac{1}{d(x,\hat{z}_i)^3} \leq C H^5 \frac{1}{d(x,\hat{z}_i)^3}. 
\end{equation}
Recalling the above definition of $\tilde{z}_i$ and $\hat{z}_i$ and the fact that 
$x_3 \in [3/4 \, H, 5/4 \, H]$, we have that 
\begin{equation}\label{eq:ddbb}
  d(x,\hat{z}_i) \geq C^{-1} H (1 + |\tilde{z}_i|^2)^{\frac{1}{2}}. 	
\end{equation}
The latter two formulas imply  
\begin{equation}\label{eq:kkcc}
 |\tilde{w}_1(x)| \leq C H^2 \sum_{\tilde{z}_i \in \Z^2}  (1 + |\tilde{z}_i|^2)^{-\frac{3}{2}} 
\leq C H^2 \int_0^\infty \frac{s}{(1+s^2)^{\frac{3}{2}}} ds \leq C H^2. 	
\end{equation}
	Concerning the gradient estimate, recalling \eqref{eq:symm-K}, the 
	fact that $K(x,y) = K(y,x)$ gives 
	\begin{equation}\label{eq:grx}
		\nabla_{\xi} K(x,y) = \nabla_{\eta} K(y,x),
	\end{equation}
	and similarly 
	\begin{equation}\label{eq:grgrx}. 
\end{equation}	
From these and 	\eqref{l:2-3}-\eqref{l:2-4} we deduce that 
$$
   |\nabla_{\xi} K(x,y)| \leq C \frac{y_3^2}{d(x,\hat{z}_i)^4}; 
  \qquad  
   |\nabla^{(2)}_{\xi,\xi} K(x,y)| \leq C \frac{y_3^2}{d(x,\hat{z}_i)^5}; 
   \qquad \quad y \in \tilde{Q}_i.
   $$
Analogously to \eqref{eq:c14}, these imply that 
$$
 \left| \int_{\tilde{Q}_i} \nabla_{\xi} K(x,y) \tilde{\theta}_1(y) dy  \right| \leq 
C H^5 \frac{1}{d(x,\hat{z}_i)^4}; \qquad \quad  \left| \int_{\tilde{Q}_i} 
\nabla^{(2)}_{\xi,\xi} K(x,y) \tilde{\theta}_1(y) dy  \right| \leq C H^5 \frac{1}{d(x,\hat{z}_i)^5}. 
$$
Since $x_3 \geq 3/4H$ and $y_3 \leq H/2$, we are avoiding the singularity of $K$ and 
we can write that 
$$
   	(\nabla \tilde{w}_1)(x) = \int_{\{0\leq y_3 \leq H/2\}} \nabla_{\xi} K(x,y) \, \th_1(y) dy; 
   	\qquad \quad 
   		(\nabla^{(2)} \tilde{w}_1)(x) = \int_{\{0\leq y_3 \leq H/2\}} \nabla^{(2)}_{\xi,\xi} K(x,y) \, \th_1(y) dy. 
$$
By the latter formulas and \eqref{eq:ddbb}, reasoning as for \eqref{eq:kkcc} we then find 
$$
 |(\nabla \tilde{w}_1)(x)| \leq C H \sum_{\tilde{z}_i \in \Z^2}  (1 + |\tilde{z}_i|^2)^{-2} 
 \leq C H \int_0^\infty \frac{s}{(1+s^2)^{2}} ds \leq C H; 
$$
$$
|(\nabla^{(2)} \tilde{w}_1)(x)| \leq C \sum_{\tilde{z}_i \in \Z^2}  (1 + |\tilde{z}_i|^2)^{-\frac{5}{2}} 
\leq C  \int_0^\infty \frac{s}{(1+s^2)^{\frac 52}} ds \leq C.  
	$$
	This concludes the proof. 	
\end{pf}

\begin{rem}\label{r:hor-deriv}
One can get estimates on higher-order derivatives of $\tilde{w}_1$ for $x_3 \in [3/4 \, H, 5/4 \, H]$ as follows. 

Consider the scaled function 
$$
  \check{w}_1(x) := \frac{1}{H^2} \hat{w}_1(H x); \qquad x = (x',x_3), \quad x'\in \R^2, x_3 \in [3/4,5/4]. 
$$
Then this is bi-harmonic, horizontally 1-periodic and by the previous proposition also uniformly bounded, together 
with its first- and second-order derivatives. Standard regularity theory implies that 
$$
  |\nabla^{(j)} \check{w}_1| \leq C_j \qquad \hbox {for } \quad  x_3 \in  [7/8, 9/8], |j| \geq 3, 
$$
which gives in turn 
\begin{equation}
	|\nabla^{(j)} \hat{w}_1| \leq C_j H^{2-j} \qquad \hbox {for } \quad  x_3 \in  [7/8 \, H, 9/8 \, H], |j| \geq 3. 
\end{equation}
\end{rem}

\section{Corrections to fit boundary data and conclusion} \label{s:pf}

We need next to adjust the above function $\tilde{w}_1$, obtained convoluting 
$K$ with $\theta_1$, see \eqref{eq:th-i},  to achieve 
the correct boundary data on $\{ x_3 = H \}$. In order to do it, we would therefore need to consider $\tilde{w}_1 - \hat{w}_1$, 
where $\hat{w}_1$ is an $H$-periodic function in $x_1, x_2$ that satisfies 
  \begin{equation}\label{eq:hat-w1}
  \begin{cases}
  \Delta^2 \hat{w}_1 = 0 & \hbox{ on }  \R^2 \times [0,H]; \\  
    \hat{w}_1 = \tilde{w}_1 & \hbox{ on } \{x_3 =  H \}; \\
    \partial_{x_3} \hat{w}_1 = \partial_{x_3} \tilde{w}_1 & \hbox{ on } \{x_3 =  H \}. 
    \end{cases}
  \end{equation}
Recall also that, by Proposition \ref{p:upper-bdry}, 
$$
  |\tilde{w}_1| \leq C \, H^2; \qquad   |\partial_{x_3} \tilde{w}_1| \leq C \, H, \qquad \quad 
  \hbox{ on } \{x_3 =  H \}. 
$$
We can solve \eqref{eq:hat-w1} via Fourier decomposition in the horizontal variables: for ${\bf k} \in \Z^2$ and ${x'} \in \R^2$,  we write 
\begin{equation}\label{eq:hatw1-Fourier}
	  \hat{w}_1({x'}, x_3) =  \sum_{{\bf k} \in \Z^2} a_{{\bf k}}(x_3) e^{2 \pi i \, \frac{\bf k}{H} \cdot {x'}}. 
\end{equation}
The function $a_{\bf k}$ satisfies the fourth-order ODE
$$
  \begin{cases}
  {a}''''_{\bf k} - 2 \left|\frac{\bf k}{H}  \right|^2 {a}''_{\bf k} + \left| \frac{\bf k}{H}  \right|^4 a_{\bf k} = 0 & \hbox{ in }  [0,H]; \\  
  a_{\bf k} = {a}'_{\bf k} = 0 & \hbox{ in } \{x_3 =  0 \}; \\
   a_{\bf k} = b_{\bf k} & \hbox{ on } \{x_3 =  H \}; \\
      a_{\bf k}' = c_{\bf k} & \hbox{ on } \{x_3 =  H \},
  \end{cases}
$$
where $b_{\bf k}$, $c_{\bf k}$ are the Fourier components of $\tilde{w}_1$ and 
$\partial_{x_3} \tilde{w}_1$ respectively on $\{x_3 =  H \}$,  defined by 
\begin{eqnarray}\label{eq:fourier-upper-boundary}
  \tilde{w}_1({x'},H) = \sum_{{\bf k} \in \Z^2} b_{{\bf k}} e^{2 \pi i \, \frac{\bf k}{H} \cdot {x'}}; 
  \qquad \qquad (\partial_{x_3} \tilde{w}_1)({x'},H) = \sum_{{\bf k} \in \Z^2} c_{{\bf k}}e^{2 \pi i \, \frac{\bf k}{H} \cdot {x'}}. 
\end{eqnarray}
From an explicit computation, see also Section C in \cite{OS} for a related one, it follows that for ${\bf k} = {\bf 0}$
 \begin{equation}\label{eq:Fourier-zero}
 a_{\bf 0} = A_{\bf 0} x_3^2 + B_{\bf 0} x_3^3; \qquad \qquad A_{\bf 0} = \frac{2 b_{\bf 0}}{H^2} - \frac{c_{\bf 0}}{2H}, \qquad 
 B_{\bf 0} = \frac{c_{\bf 0}}{2H^2} - \frac{b_{\bf 0}}{H^3}, 
 \end{equation}
while for ${\bf k} \neq {\bf 0}$
\begin{equation}\label{eq:ak-zeta}
  a_{\bf k}(x_3) = x_3 (B_{\textbf{k}} \sinh (|{\bf k} |/H \, x_3)-A_{\textbf{k}} |\textbf{k}|/H \cosh (|{\bf k} |/H \, x_3))+A_{\textbf{k}} \sinh (|{\bf k} |/H \, x_3), 
\end{equation}
where 
$$
  A_{\bf k}  = \frac{2 ((b_{\bf k} -c_{\bf k}  H) \sinh (|{\bf k} |)+b_{\bf k}  |{\bf k} |  \cosh (|{\bf k} |))}{\cosh (2 |{\bf k} |)-2 |{\bf k} |^2-1}; 
 $$ 
 $$
 B_{\bf k}  = 2 \frac{\sinh (|{\bf k} |) \left(b_{\bf k}  |{\bf k} |^2/H+c_{\bf k} \right) - c_{\bf k}  |{\bf k} |  \cosh (|{\bf k} | )}{\cosh (2 |{\bf k} |)-2 |{\bf k} |^2-1}. 
$$
Due to Remark \ref{r:hor-deriv}, the Fourier components of $\tilde{w}_1$ and $\partial_{x_3} \tilde{w}_1$ for $x_3 =  H$ 
decay fast in ${\bf k}$: more precisely, for any integer $\ell$ there exists $C_\ell > 0$, independent of $H$ large, 
such that 
\begin{equation}\label{eq:fourier-bounds}
  \sum_{{\bf k} \in \Z^2} \frac{|b_{{\bf k}}/H^2|^2 + |c_{{\bf k}}/H|^2}{1 + |{\bf k}|^\ell} \leq C_\ell.  
\end{equation}
This will imply 
uniform convergence of the series (in ${\bf k}$) of $\hat{w}_1$ for $x_3 \in [0,1]$ and of that of $\tilde{w}_1 - \hat{w}_1$ for $x_3 \in [H-1,H]$. We indeed prove the 
following result.

\begin{pro}\label{p:bihar-bd}
	Let $\hat{w}_1$ be as in \eqref{eq:hat-w1}. 
	There exist fixed constants $H_0$ and $C_0$ such that for all $H \geq H_0$ 
	$$
	  |\hat{w}_1| \leq C_0 \quad \hbox{ in } \{ 0 \leq x_3 \leq 1 \}; \qquad \quad 
	   |\tilde{w}_1 - \hat{w}_1| \leq C_0 \quad \hbox{ in }\{ H-1 \leq x_3 \leq H \}. 
	$$
\end{pro}

\begin{pf}
We first estimate the zero-mode $a_{\bf 0}$ in \eqref{eq:Fourier-zero}. For $x_3 \in [0,1]$, $a_{\bf 0}(x_3)$ is clearly bounded since 
$A_{\bf 0}$ and $B_{\bf 0}$ are, due to Proposition \ref{p:upper-bdry}. 

Let us now pass to higher-order modes. We first  notice that, since the denominators in $A_{\bf k}$, $B_{\bf k}$ are bounded below by $C^{-1} e^{2 |{\bf k}| }$, 
then 
\begin{equation}\label{eq:bds-Ak-Bk}
	|A_{\bf k} | \leq C (|b_{\bf k}| + H |c_{\bf k}|) |{\bf k}|  e^{- |{\bf k}| }; \qquad 
	\quad |B_{\bf k} | \leq C (|b_{\bf k}|/H + |c_{\bf k}|) |{\bf k}|^2  e^{- |{\bf k}| }. 
\end{equation}
With a change of variables, we can write 
\begin{equation}\label{eq:akz}
	  a_{\bf k}(x_3) = \left( A_{\bf k}  + \frac{H}{|{\bf k}|} B_{{\bf k}} \, z \right) \sinh z - z A_{{\bf k}} \cosh z; 
	  \qquad z = \frac{|{\bf k}|}{H} x_3. 
\end{equation}
This function vanishes at $z = 0$ and has also vanishing first-order derivative at $z= 0$, while 
its second-order derivative is given by 
$$
  2 \frac{H}{|{\bf k}|} B_{{\bf k}} \cosh z + \frac{H}{|{\bf k}|} B_{{\bf k}} z \sinh z - A_{{\bf k}} (\sinh z + z \cosh z). 
$$
By a Taylor expansion in integral form, one finds that 
$$
 |a_{\bf k}(x_3)| \leq C z^2 e^z \left( |A_{{\bf k}}| + \frac{H}{|{\bf k}|} |B_{\bf k}| \right); \qquad z = \frac{|{\bf k}|}{H} x_3, 
$$
which implies 
$$
|a_{\bf k}(x_3)| \leq C \left( \frac{|{\bf k}|}{H} \right)^2 e^{\frac{|{\bf k}|}{H}} \left( |A_{{\bf k}}| + \frac{H}{|{\bf k}|} |B_{\bf k}| \right); 
\qquad x_3 \in [0,1]. 
$$
From \eqref{eq:bds-Ak-Bk}, it follows that 
$$
  |a_{\bf k}(x_3)| \leq C \left( \frac{|{\bf k}|}{H} \right)^2 e^{\frac{|{\bf k}|}{H}} (|b_{\bf k}| + H |c_{\bf k}|) |{\bf k}|  e^{- |{\bf k}| }
  \leq C |{\bf k}|^3 e^{- \frac{1}{2} |{\bf k}| } \left(  \frac{|b_{\bf k}|}{H^2} +  \frac{|c_{\bf k}|}{H} \right); 
  \qquad x_3 \in [0,1], 
$$
%
%
so by \eqref{eq:fourier-bounds} the series in \eqref{eq:hatw1-Fourier}   converges arbitrarily fast for  $x_3 \in [0,1]$. 
This proves the first statement in the proposition. 

\

We turn next to the second assertion. It will be sufficient to 
control the second derivative in $x_3$ of $\hat{w}_1 $. 
Using \eqref{eq:akz}, we have that 
\begin{eqnarray*}
	\frac{d^2}{d x_3^2} a_{{\bf k}}(x_3) & = & \left( \frac{dz}{d x_3} \right)^2 \frac{d^2}{d z^2} a_{{\bf k}}(x_3) 
	\\ & = & \left( \frac{|{\bf k}|}{H} \right)^2 \left[ 2 \frac{H}{|{\bf k}|} B_{{\bf k}} \cosh z - A_{{\bf k}} \sinh z  
	+ z \left( \frac{H}{|{\bf k}|} B_{{\bf k}}  \sinh z - A_{{\bf k}} \cosh z \right)
	\right]. 
\end{eqnarray*}
To evaluate this quantity for $x_3 \in [H-1,H]$, we write 
\begin{eqnarray*}
 & & 2 \frac{H}{|{\bf k}|} B_{{\bf k}} \cosh z - A_{{\bf k}} \sinh z  
 + z \left( \frac{H}{|{\bf k}|} B_{{\bf k}}  \sinh z - A_{{\bf k}} \cosh z \right) \\ & = & 
 \left[ \left(  2 \frac{H}{|{\bf k}|} B_{{\bf k}} - A_{\bf k}  \right) 
 + z \left(  \frac{H}{|{\bf k}|} B_{{\bf k}} - A_{\bf k}  \right)  \right] \sinh z  + 
 \left(  2 \frac{H}{|{\bf k}|} B_{{\bf k}} - A_{\bf k}  \right)  e^{-z} - z \left(   \frac{H}{|{\bf k}|} B_{{\bf k}} - A_{\bf k}  \right)  e^{-z}. 
\end{eqnarray*}
For $x_3 \in [H-1,H]$, by \eqref{eq:bds-Ak-Bk} the latter expression can be bounded by 
$$
 C z e^{z} \left(   \frac{H}{|{\bf k}|} |B_{{\bf k}}|  +  |A_{\bf k}|  \right) \leq C |{\bf k}|  e^{|{\bf k}|}  
  (|b_{\bf k}| + H |c_{\bf k}|) |{\bf k}|  e^{- |{\bf k}| } \leq C |{\bf k}|^2  (|b_{\bf k}| + H |c_{\bf k}|).  
$$
 The above formula for $\frac{d^2}{d x_3^2} a_{{\bf k}}(x_3) $ then implies 
 $$
   \left| \frac{d^2}{d x_3^2} a_{{\bf k}}(x_3)  \right| \leq C \frac{|{\bf k}|^4}{H^4} \frac{d^2}{d x_3^2} a_{{\bf k}}(x_3); 
   \qquad x_3 \in [H-1,H]. 
 $$
By \eqref{eq:fourier-bounds} this function is uniformly bounded, giving the conclusion 
by  the last estimate in \eqref{eq:estx3L} and by the fact that $\tilde{w}_1 - \hat{w}_1$ 
vanishes and has zero normal derivative on the plane $\{ x_3 = H\}$. 
\end{pf}

We can now prove our main result. 

%
%

\

\begin{pfn} {\sc of Theorem \ref{t:main}}. 
	Let $\th_1$, $\th_2$ be as in \eqref{eq:th-i}, and let $w_i$, $i = 1, 2$, be the solution of 
$$
\begin{cases}
	\D^2 w_i = \Delta_h \th_i & \hbox{ in } \{ 0 \leq x_3 \leq H \}; \\ 
	w_i = 0 & \hbox{ on } \{ x_3 =  0 \} \cup  \{ x_3 =  H \}; \\ 
	\partial_{x_3} w_i = 0 & \hbox{ on } \{ x_3 =  0 \} \cup  \{ x_3 =  H \}. 
\end{cases}
$$
If $\tilde{w}_1$ is as in \eqref{eq:repr} and $\hat{w}_1$ as in \eqref{eq:hat-w1}, 
then by uniqueness 	$w_1 = \tilde{w}_1 - \hat{w}_1$. Proposition \ref{p:w1} (recalling 
that $\th = \th_1$ for $x_3 \leq H/2$) and the first estimate in Proposition \ref{p:bihar-bd} imply that 
$\langle \th \, w_1 \rangle$ is uniformly bounded for $x_3 \in (0,1]$. 

On the other hand, we could apply the second estimate in Proposition \ref{p:bihar-bd} 
to the function $w_2(x',H-x_3)$ and again use the $L^\infty$-bound on $\theta$ 
to show that also $\langle \th \, w_2 \rangle$ is uniformly bounded for $x_3 \in (0,1]$. 
The conclusion  follows then from formula \eqref{eq:bd-Nu-OS}. 
\end{pfn}

	\section{Appendix} \label{s:appendix}

	In this appendix we collect some remarks to help the reader understand the proof 
	of Proposition \ref{p:CM-mod} and to point out the elementary observations that allow us to adapt the original one for Theorem 33 in \cite{CM}.
	In fact this modification is already suggested by the  proof of that result.  
	We collect first some facts from the 
	 book \cite{St0}.

	 \begin{lem}[\cite{St0}, page 62] \label{l:app-1}
	 	Let $P_t$ denote the Poisson kernel, see \eqref{eq:Poisson}. Then its Fourier transform 
	 	has the expression 
	 	$$
	 	\hat{P}_t(\xi) = e^{- 2 \pi t |\xi|}, \qquad \xi \in \R^2. 
	 	$$
	 	In particular, $\hat{P}_t$ is exponentially decaying at infinity and is in the Schwartz space at infinity. 
	 \end{lem}

	It is useful to recall the following definition.

	\begin{df} \label{df:carleson}
		Let $\mu$ be a measure on $\R^3_+$. Let us consider the {\em Carleson box} 
		$$
		 T(x^0,h) := \left\{ (x,h) \in \R^2 \times \R \; | \; |x-x^0| < h, 0 \leq t < h \right\}.
		$$
		We say that $\mu$ is a {\em Carleson measure} if there is an 
		absolute constant $C_1$ such that   for all $h > 0$ 
		$$
		 \mu(T(x^0,h)) \leq C_1 h^2. 
		$$
	\end{df}

			\begin{figure}[!htb]
				\begin{center}
					\minipage{0.80\textwidth}
					\includegraphics[width=\linewidth]{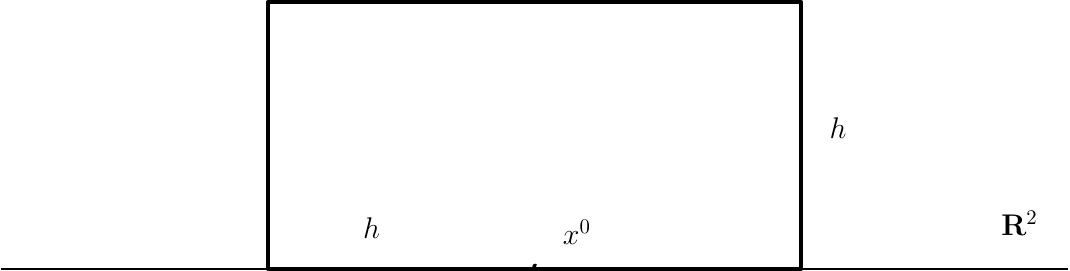}
					\caption{\small A Carleson box}\label{fig:box}
					\endminipage\hfill
				\end{center} 
			\end{figure}

	Recalling the definition of 		non-tangential maximal function in \eqref{eq:NTM}, we have the following result.

	\begin{lem}[\cite{St0}, page 236] \label{l:stein}
		Let $F(x,t)$ be any measurable function on $\R^3_+$, and let $\mu$ be a 
		Carleson measure. Then, for $0 < p < + \infty$  
		$$
		\int_{\R^3_+} |F(x,t)|^p d \mu \leq C(p) \, C_1 \int_{\R^2} |N  F|^p dx, 
		$$ 
		where $C_1$ is the constant in  Definition \ref{df:carleson} and and $N F$ is the 
			non-tangential maximal function of $F$. 
	\end{lem}

The case $p = 2$ is  employed in \cite{CM}. We will  call {\em bump functions} 
approximate identities  of the type 
$$
 \varphi_t(x) = \frac{1}{t^n} \varphi (x/t).
$$
Even though Proposition \ref{p:CM-mod} is stated for vectorial kernels, it will be enough to 
consider scalar ones. We will  employ as bump functions  $t \nabla_h P_t$, where $\nabla_h P_t$ stands for  the horizontal derivatives 
of the Poisson kernel $P_t$ on $\R^2$. The Fourier transform of this function when $t = 1$ 	is 
$$
 \xi_j \hat{P}(\xi) = \xi_j e^{- 2 \pi |\xi|},
$$
see e.g. \cite[page 61]{St0}.
We have the following estimates, for all multi-indices $\alpha$
\begin{equation}
	\begin{cases} \label{eq:111}
		| \partial^\alpha_\xi \hat{P}(\xi)| \leq C(\alpha) e^{- C(\alpha)^{-1} |\xi|}  & \hbox{ for } |\xi| \geq 1; \\ 
		|\xi|^\alpha |\partial^\alpha_\xi  \hat{P}(\xi)| \leq C(\alpha) |\xi| & \hbox{ for } |\xi| \leq  1. 
	\end{cases}
\end{equation}	
Thus the above inequalities imply the decay assumption in Theorem 33 of \cite{CM}, reported in Remark \ref{r:CM-gen}. 

The proof of Theorem 33 in \cite{CM} has a key step, namely  Lemma 5 on page 151. Here one establishes the theorem under the provisional assumption that the bump functions have 
compactly supported Fourier transform,  with the origin not lying in the support of the Fourier transform of any of the bump functions, and also the origin not lying in the algebraic sum of the supports of the Fourier transform of the bump functions. 

In the general case  where  assumptions \eqref{eq:111} hold true, the arguments on pages 152-153 in \cite{CM} 
reduce matters to Lemma 5 there. We now make some remarks on the proof of Lemma 5, which is already implicit in the proof displayed in \cite{CM}, but that we now highlight. 
We collect real variable facts used in Lemma 5 which are standard.

\begin{lem} \label{l:l3-app}
	For two functions $\varphi, \psi \in L^1(\R^2)$ one has 
	\begin{description}
		\item[(i)] $\widehat{\varphi * \psi}(\xi) = \hat{\varphi}(\xi) \hat{\psi}(\xi)$; 
		\item[(ii)] $\widehat{\varphi \psi}(\xi) = (\hat{\varphi}*\hat{\psi})(\xi)$; 
		\item[(iii)] the support of $\hat{\varphi}*\hat{\psi}$ is contained in the algebraic sum of 
		\textup{supp} $\hat{\varphi}$ and  \textup{supp} $\hat{\psi}$. 
	\end{description}
\end{lem}

In \cite{CM}, one considers integrals of the type 
$$
  \int_0^\infty (f * \varphi_t)(x) (b * \psi_t)(x) m(t) \frac{dt}{t}, 
$$
namely 
$$
\int_0^\infty (F(\cdot,t) * \varphi_t)(x) (b * \psi_t)(x) \frac{dt}{t},   
$$
with 
 \begin{equation}\label{eq:F-CM}
   {F}(x', t) = f(x' ) m(t).  	
 \end{equation}
In our case, for Proposition \ref{p:CM-mod}, we need to estimate integrals of the type   
$$
\int_0^\infty (F(\cdot,t)  * \varphi_t)(x) (b * \psi_t)(x) \frac{dt}{t},   
$$
with $\varphi_t  = \psi_t$  and a \underline{more general dependence of $F$ on $t$}, compared to \eqref{eq:F-CM}.

The latter integral is a paraproduct, see \cite{CM}, Appendix I. To estimate the $L^2(\R^2)$-norm of the  above expression, in  
\cite{CM} one proceeds via duality and consider, for $h$ of class $L^2(\R^2)$, the operator 
\begin{equation}\label{eq:app-3}
	h \quad \longmapsto \quad \int_{\R^2} h(x) \left[  \int_0^\infty (F(\cdot, t) * \varphi_t)(x) (b * \psi_t)(x) m(t) \frac{dt}{t} \right] dx.  
\end{equation}
The main idea of Lemma 5 in \cite{CM} (see page 151 there), under the above-mentioned assumptions on the 
supports of $\hat{\varphi}$ and $\hat{\psi}$, is that one can write, for $\delta > 0$ 
$$
[(F(\cdot, t) * \varphi_t) (\cdot) (b * \psi_{\delta t})(\cdot)] m(t) \qquad 
\hbox{ as } \qquad 
\psi_{3,t} * [(F(\cdot, t) * \varphi_t) (\cdot) (b * \psi_{\delta t})(\cdot)] m(t). 
$$
Here $\psi_{3,t}(x)  = \frac{1}{t^2} \psi_3(x/t)$ 
is a bump function 
for which $\hat{\psi}_3$ is compactly supported, smooth, such that 
$\hat{\psi}_3 \equiv 1$ on $\textup{supp } \hat{\varphi}(\cdot) + \textup{supp } \hat{\psi}(\delta \cdot)$ and 
with  $\textup{supp } \hat{\psi}_3$ contained either in a 
ball or in a proper  annulus around the origin. 
When applying then $L^2$-duality in \eqref{eq:app-3}, one can then pass the convolution with $\psi_{3,t} $ 
to $h$ and use  Littlewood-Paley estimates. 

To explain these steps in  more detail,  consider the Fourier transform in $x$ of 
\begin{equation}\label{eq:app-4}
	(F(\cdot, t)*\varphi_t)(x) (b * \psi_{\delta t})(x) m(t). 
\end{equation}
By Lemma \ref{l:l3-app}, this is given by 
\begin{equation}\label{eq:app-5}
	[\hat{F}(\xi,t) \widehat{\varphi_t}(\xi)] * [\hat{b}(\xi) \widehat{\psi_t}(\xi)] m(t).  
\end{equation}
The point is now that the support of $\hat{F}(\xi,t) \widehat{\varphi_t}(\xi)$ lies in the support of  
 $\widehat{\varphi_t}$ and the support of $\hat{b}(\xi) \widehat{\psi_t}(\xi)$ lies in the 
support of $\widehat{\psi_t}$, and for this to hold it is irrelevant that $F(x,t)$ depends on $t$ or that there is 
the extra factor  $m(t)$. Thus the support of \eqref{eq:app-5} lies in the algebraic sum of the supports 
of $\widehat{\varphi_t}$ and  $\widehat{\psi_t}$.


By the above discussion, the dual pairing of $h$ and the expression in \eqref{eq:app-4} becomes 
$$
 	\int_{\R^2} h(x)   \int_0^\infty  \psi_{3,t} *[   (F(\cdot, t) * \varphi_t)  (b * \psi_{\delta t})](x) m(t) \frac{dt}{t}  dx, 
$$
and after switching convolutions it can be written as 
$$
\int_{\R^2} \int_0^\infty (\psi_{3,t} *h)(x)     [ (F(\cdot, t) * \varphi_t)  (b * \psi_{\delta t})](x) m(t) \frac{dt}{t}  dx.  
$$
Using the fact that $m \in L^\infty$ and the Cauchy-Schwarz inequality, this quantity can be 
bounded by 
\begin{eqnarray} \label{eq:61}
	\left( \int_{\R^2} \int_0^\infty |(\psi_{3,t}*h)(x)|^2 \frac{dt}{t} dx \right)^{\frac 12 } 
	 \times  \left( \int_{\R^2} \int_0^\infty  |(F(\cdot, t)*\varphi_t)(x)|^2 |(b * \psi_{\delta t})(x)|^2 \frac{dt}{t} dx \right)^{\frac 12}.
\end{eqnarray}
We apply Lemma 1 in \cite{CM} to the first term in the latter formula to conclude that it is bounded by a fixed constant times $\|h\|_{L^2(\R^2)}$.
To the second term we apply Lemma \ref{l:stein}. 
Recall that we set 
$$
 (F(\cdot, t) * \varphi_t)(x) = \check{F}(x,t): 
$$
from Theorem \ref{t:FS} we see that  $|(b*\psi_t)(x)|^2 \frac{dt}{t} dx$ is a 
Carleson measure, so the second term in \eqref{eq:61}, using Lemma \ref{l:stein}, can be bounded by 
$$
  \left( \int_{ \R^2} (N  \check{F})^2(x) dx \right)^{\frac 12}, 
$$
as desired. For general kernels as in Proposition \ref{p:CM-mod}, see the comments before Lemma \ref{l:l3-app}, we 
 follow the arguments on pages 152-153 in \cite{CM}.


\end{document}